\documentclass{amsart}
\begin{document}

\vfuzz2pt 
\hfuzz2pt 
\newtheorem{thm}{Theorem}[section]
\newtheorem{corollary}[thm]{Corollary}
\newtheorem{lemma}[thm]{Lemma}
\newtheorem{proposition}[thm]{Proposition}
\newtheorem{defn}[thm]{Definition}
\newtheorem{remark}[thm]{Remark}
\newtheorem{example}[thm]{Example}
\newtheorem{fact}[thm]{Fact}
\
\newcommand{\norm}[1]{\left\Vert#1\right\Vert}
\newcommand{\abs}[1]{\left\vert#1\right\vert}
\newcommand{\set}[1]{\left\{#1\right\}}
\newcommand{\Real}{\mathbb R}
\newcommand{\eps}{\varepsilon}
\newcommand{\To}{\longrightarrow}
\newcommand{\BX}{\mathbf{B}(X)}
\newcommand{\A}{\mathcal{A}}
\newcommand{\onabla}{\overline{\nabla}}
\newcommand{\hnabla}{\hat{\nabla}}


\def\proof{\medskip Proof.\ }
\font\lasek=lasy10 \chardef\kwadrat="32 
\def\kwadracik{{\lasek\kwadrat}}
\def\koniec{\hfill\lower 2pt\hbox{\kwadracik}\medskip}

\newcommand*{\C}{\mathbf{C}}
\newcommand*{\R}{\mathbf{R}}
\newcommand*{\Z}{\mathbf {Z}}

\def\sb{f:M\longrightarrow \C ^n}
\def\det{\hbox{\rm det}\, }
\def\detc{\hbox{\rm det }_{\C}}
\def\i{\hbox{\rm i}}
\def\tr{\hbox{\rm tr}\, }
\def\rk{\hbox{\rm rk}\,}
\def\vol{\hbox{\rm vol}\,}
\def\Im {\hbox{\rm Im}\, }
\def\Re{\hbox{\rm Re}\, }
\def\interior{\hbox{\rm int}\, }
\def\e{\hbox{\rm e}}
\def\pu{\partial _u}
\def\pv{\partial _v}
\def\pui{\partial _{u_i}}
\def\puj{\partial _{u_j}}
\def\puk{\partial {u_k}}
\def\div{\hbox{\rm div}\,}
\def\Ric{\hbox{\rm Ric}\,}
\def\r#1{(\ref{#1})}
\def\ker{\hbox{\rm ker}\,}
\def\im{\hbox{\rm im}\, }
\def\I{\hbox{\rm I}\,}
\def\id{\hbox{\rm id}\,}
\def\exp{\hbox{{\rm exp}^{\tilde\nabla}}\.}
\def\cka{{\mathcal C}^{k,a}}
\def\ckplusja{{\mathcal C}^{k+1,a}}
\def\cja{{\mathcal C}^{1,a}}
\def\cda{{\mathcal C}^{2,a}}
\def\cta{{\mathcal C}^{3,a}}
\def\c0a{{\mathcal C}^{0,a}}
\def\f0{{\mathcal F}^{0}}
\def\fnj{{\mathcal F}^{n-1}}
\def\fn{{\mathcal F}^{n}}
\def\fnd{{\mathcal F}^{n-2}}
\def\Hn{{\mathcal H}^n}
\def\Hnj{{\mathcal H}^{n-1}}
\def\emb{\mathcal C^{\infty}_{emb}(M,N)}
\def\M{\mathcal M}
\def\Ef{\mathcal E _f}
\def\Eg{\mathcal E _g}
\def\Nf{\mathcal N _f}
\def\Ng{\mathcal N _g}
\def\Tf{\mathcal T _f}
\def\Tg{\mathcal T _g}
\def\diff{{\mathcal Diff}^{\infty}(M)}
\def\embM{\mathcal C^{\infty}_{emb}(M,M)}
\def\U1f{{\mathcal U}^1 _f}
\def\Uf{{\mathcal U} _f}
\def\Ug{{\mathcal U} _g}
\def\[f]{{\mathcal U}^1 _{[f]}}
\def\hnu{\hat\nu}
\def\gnu{\nu_g}
\title{A sectional curvature for statistical structures}
\author{Barbara Opozda}



\thanks{The research supported by the NCN grant K/PBO/000302 and a grant of the TU in Berlin } \maketitle

\address{Instytut Matematyki UJ, ul. \L ojasiewicza  6, 30-348
Cracow, Poland}

\email{Barbara.Opozda@im.uj.edu.pl}


 \begin{abstract}{A new type of sectional curvature is introduced.
 The notion is purely algebraic and can be located
   in linear algebra as well as in differential
 geometry.}
\end{abstract}

\maketitle

\section{Introduction}

Sectional curvature is one of the most important concepts in
differential geometry. Nevertheless, it is attributed to Riemannian
or pseudo-Riemannian geometry only. The curvature tensor field is
defined for any connection but to define a sectional curvature,
which assigns to a vector plane of a tangent space a number, seems
to need a scalar product. Moreover, the metric and the connection
must be related in a good manner. For instance, in the classical
affine differential geometry one has a metric tensor field and the
so called induced connection, but the curvature tensor of type
$(0,4)$ constructed by these objects does not have enough
symmetries.  The tensor satisfies appropriate symmetry conditions
for affine spheres but it leads to  trivial cases, namely to spaces
of constant sectional curvature. The problem can be  solved by
adding to the curvature tensor the curvature tensor for the dual
connection.  This idea is discussed in \cite{O} for statistical
structures on abstract manifolds, that is,  on  manifolds (not
necessarily immersed in any standard space) endowed with a matric
tensor field $g$ and a torsion-free affine connection $\nabla$ for
which $\nabla g$ as a 3-covariant tensor is symmetric.

A statistical structure is also called a Codazzi structure, see e.g.
\cite{NSi}, \cite{NS}. We use the name  ''statistical structure''
following \cite{MTA} or \cite{L}. The name  ''Codazzi structure''
may refer to all situations, where we have any tensor field whose
covariant derivative is totally symmetric.

The  geometry of affine hypersurfaces in the standard affine space
$\R^n$ or, more generally, the geometries of the second fundamental
form, including the theory of Lagrangian submanifolds in complex
space forms,  are natural sources of statistical structures.
However, the fact that the structures are induced by the simple
structures on the ambient spaces imposes strong conditions on the
induced statistical structure. For instance,  for affine
hypersurfaces, it it necessary that the dual connection is
projectively flat.

It turns out that for statistical structures one can define few
sectional curvatures. In \cite{O}  we studied the sectional
$\nabla$-curvature, that is, a sectional curvature  determined by a
metric tensor and a connection $\nabla$.
 In
this paper we propose  another type of sectional curvature. Its idea
is purely algebraic. This sectional curvature can be defined on any
 vector space endowed with a scalar product and a symmetric cubic
form. Then it can be transfered to statistical structures on
manifolds. In this paper we provide some basic information on this
sectional curvature and we give exemplary theorems concerning this
notion.

\section{Statistical structures}
One can define a statistical structure on a  manifold $M$  in three
equivalent ways. First of all $M$  must have a Riemannian structure
defined by a metric tensor field $g$. Throughout  the paper we
assume that $g$ is positive definite, although $g$ can be also
indefnite. A statistical structure can be defined as a pair $(g,K)$
on a manifold $M$, where $g$ is a Riemannian metric tensor field and
$K$ is a symmetric $(1,2)$-tensor field which is also symmetric
relative to $g$, that is, the cubic form
\begin{equation}\label{K_C}
C(X,Y,Z)=g(X,K(Y,Z))
\end{equation}
 is symmetric relative to $X,Y$.  It is clear that any
symmetric cubic form $C$ on a Riemannian manifold $(M, g)$ defines
by (\ref{K_C}) a $(1,2)$-tensor field $K$ having the symmetry
properties as above. Another equivalent definition says that a
statistical structure is a pair $(g,\nabla)$, where $\nabla$ is a
torsion-free affine connection on $M$ and $\nabla g$ as a
$(0,3)$-tensor field on $M$ is symmetric in all arguments. Let us
fix  that for a tensor field $s$ and a connection $\nabla$ the
notation $\nabla s (X,...)$ stands for $(\nabla _Xs)(...)$. The
affine connection $\nabla$ from the last definition equals to
$\hat\nabla +K$, where $\hat\nabla$ is the Levi-Civita connection
for $g$ and $K$ is the difference tensor. Since  $\nabla g(X,
Y,Z))=-2g(K(X,Y),Z)$, we obtain a statistical structure $(g,K)$ from
$(g,\nabla)$. We shall call $\nabla$ a statistical connection. A
manifold equipped with a statistical  structure will be called a
statistical manifold.

For any connection $\nabla$ on a Riemannian manifold $(M,g)$ one
defines its conjugate connection $\onabla$ (relative to $g$) as
follows
\begin{equation}
g(\nabla _XY,Z)+g(Y,\onabla _XZ)=Xg(Y,Z)
\end{equation}
for any vector fields $X,Y,Z$ on $M$. The connections $\nabla$ and
$\overline\nabla$ are simultaneously torsion-free.  It is also known
that if $(g,\nabla)$ is a statistical structure then so is
$(g,\onabla)$. Moreover, if $(g,\nabla)$ is trace-free then so is
$(g,\onabla)$, see e.g. \cite{NS}. Recall that a trace-free
statistical structure is such a structure for which $\tr _g(\nabla
g)(X,\cdot, \cdot)=0$ for every $X$ or equivalently $\tr _gK=0$, or
equivalently $\tr K_X=0$  for every $X$, where $K_XY=K(X,Y)$. Note
that a statistical structure is trace-free if and only if $\nabla
\nu_g=0$, where $\nu_g$ is the volume form determined by $g$. If $R$
is the curvature tensor for $\nabla$ and $\overline R$ is the
curvature tensor for $\onabla$ then we have, \cite{NS},
\begin{equation}\label{R_and_oR}
g(R(X,Y)Z,W)=-g(\overline{R}(X,Y)W,Z)
\end{equation}
for every $X,Y,Z,W$. In particular,  $R=0$ if and only if $\overline
R=0$. If $K$ is the difference tensor between $\nabla$ and
$\hnabla$, that is,
\begin{equation}
\nabla _XY=\hnabla _XY+K_XY,
\end{equation}
then
\begin{equation}
\onabla_XY=\hnabla _XY-K_XY.
\end{equation}
 It is also known that
\begin{equation}\label{from_Nomizu_Sasaki}
R(X,Y)=\hat R(X,Y) +(\hnabla_XK)_Y-(\hnabla_YK)_X+[K_X,K_Y].
\end{equation}
Writing the same equality for $\onabla$ and adding both equalities
we get
\begin{equation}\label{R+oR}
R(X,Y)+\overline R(X,Y) =2\hat R(X,Y) +2[K_X,K_Y].
\end{equation}
 The
following lemma follows from  formulas (\ref{R_and_oR}),
(\ref{from_Nomizu_Sasaki}) and (\ref{R+oR}).

\begin{lemma}\label{przeniesiony_lemat}
Let $(g,K) $ be a statistical structure. The following conditions
are equivalent:
\newline
{\rm 1)} $R=\overline R$,
\newline
{\rm 2)} $\hnabla K$ is symmetric,
\newline
{\rm 3)} $g(R(X,Y)Z,W)$ is skew-symmetric for $Z,W$.
\end{lemma}

A statistical structure is called  Hessian  if the connection
$\nabla$ is flat, that is, $R=0$. In this case, by (\ref{R+oR}), we
have
\begin{equation}
\hat R=-[K,K].
\end{equation}
For a statistical structure one  defines  the  vector field $E$ by
\begin{equation}
E=\tr_g K.
\end{equation}
If $e_1,...,e_n$ is an orthonormal basis of $\mathcal V$ then
\begin{equation}\label{trK}
E=(\tr K_{e_1})e_1 +...+(\tr K_{e_n})e_n.
\end{equation}
 For more information on dual connections, affine differential
geometry and statistical structures we refer to \cite{NSi},
\cite{LSZ}, \cite{NS}, \cite{L}, \cite{MTA}, \cite{Sh}, \cite{O}.

\section{The sectional $K$-curvature}
First we shall give an  algebraic setting of the sectional
$K$-curvature. Let $\mathcal V$ be a vector space with a positive
definite scalar product $g$. Let $K$ be a symmetric tensor field of
type $(1,2)$ on $\mathcal V$ and symmetric relative to $g$.
 Hence $K_X$ is a tensor of type $(1,1)$ symmetric relative
to $g$. In particular, it is diagonalizable. $K$ defines a symmetric
cubic form $C$ given by (\ref{K_C}).

The tensor field $K$ determines a $(1,3)$-tensor  $[K,K]$ given by
$$[K,K](X,Y)Z:=[K_X,K_Y]Z=K_XK_YZ-K_YK_XZ.$$ This is a
curvature-like tensor, that is, it satisfies the following
conditions
\begin{eqnarray*}
&&  [K,K](X,Y)=-[K,K](Y,X)\\
&&
 [K,K](X,Y)Z+[K,K](Y,Z)X +[K,K] (Z,X)Y=0\\
 &&
 g([K,K](X,Y)Z,W)=-g([K,K](X,Y)W,Z)
\end{eqnarray*}
for every vectors $X,Y,Z,W\in\mathcal V$. It follows that  we can
define the sectional $K$-curvature
by a vector plane $\pi$ in $\mathcal V$ as follows.
Take an orthonormal basis $X,Y$ of $\pi$ and set
\begin{equation}
 k(\pi) = g([K,K](X,Y)Y,X).
\end{equation}
The number $k(\pi)$ is independent of the choice of an orthonormal
basis $X,Y$. The sectional $K$-curvature by a plane spanned by
vectors $X,Y$ will be denoted by $k(X\wedge Y$).

On a $2$-dimensional vector  space $\mathcal V$ we have
$[K,K](X,Y)Z=k(\mathcal V)[g(Y,Z)X-g(X,Z)Y]$ for all vectors
$X,Y,Z\in \mathcal V$. If  the dimension of $\mathcal V$ is
arbitrary and the sectional $K$-curvature is  equal to some constant
number $A$ for all vector planes in $\mathcal V$ then we have
\begin{equation}\label{constant_curvature_1}
 [K,K](X,Y)Z=A [g(Y,Z)X-g(X,Z)Y]
\end{equation}
for every $X,Y,Z\in \mathcal V$. The condition
(\ref{constant_curvature_1}) can be written equivalently as
\begin{equation}\label{constant_curvature_2}
\begin{array}{rcl}
&&\ \ g(K(X,W), K(Y,Z))-g(K(Y,W),K(X,Z))\\
&&\ \ \ \ \ \ \ \ \ \ \ \ \ \ \ \ \ \ \ \ \ \ \ \ \ =
A[g(X,W)g(Y,Z)-g(Y,W)g(X,Z)]
\end{array}
\end{equation}
for every $X,Y,Z,W\in \mathcal V$.

The sectional $K$-curvature can be now introduced on a statistical
manifold $(M,g,K)$ in the above manner on each tangent space. In
general, Schur's lemma does not hold. It follows from the fact that
the curvature tensor $[K,K]$, in general, does not satisfy any
second Bianchi identity. The following identity can be regarded as
the second Bianchi identity for the curvature tensor $R+\overline
R$, see \cite{O},
\begin{lemma}\label{IIBianchi}
For any statistical structure $(g,\nabla)$ we have
$$\Xi _{U,X,Y}(\hnabla _U(R+\overline R))(X,Y)=
\Xi _{U,X,Y}( K_U\cdot(\overline R- R))(X,Y),$$ where $\Xi$ stands
for the cyclic permutation sum.
\end{lemma}
 Using Lemmas (\ref{przeniesiony_lemat} and (\ref{IIBianchi}) one
easily gets the following analogue of Schur's lemma
\begin{proposition}
Let $(g,K)$ be a statistical structure on a connected manifold $M$
whose dimension is greater than 2. If the  (1,3)-tensor field
$\hat\nabla K$ is symmetric  and the sectional $K$-curvature depends
only on a point of $M$ then  the sectional $K$-curvature is constant
on $M$.
\end{proposition}

\begin{example}\label{example_k=lambda_2over4}
{\rm Let $e_1,...,e_n$ be an orthonormal frame  of $\mathcal V$.
Define a $(1,2)$-tensor $K$ on $\mathcal V$ as follows
\begin{equation}
 K(e_1,e_1)=\lambda e_1, \ \ \ \ K(e_1,e_i)=
 \frac{\lambda}{2} e_i, \ \ \ \ K(e_i,e_i)=\frac{\lambda}{2} e_1,\ \ \ \
 K(e_i,e_j)=0
\end{equation}
for $i,j\ge 2$, $i\ne j$. By a straightforward computation one can
check that the sectional $K$-curvature  is constant on $\mathcal V$
and equals to $\lambda^2/4$. Observe that in this case the cubic
form $C$ vanishes on the $(n-1)$-dimensional hyperplane spanned by
$e_2,..., e_n$ and the sectional curvature is positive for all
sections.
 }
\end{example}
\begin{example}{\rm In \cite{Ch} B-Y Chen studied Lagrangian
$H$-umbilical submanifolds. That study leads to the following
examples of statistical structures. Let $g$ be a scalar product on a
vector space $\mathcal V$. In some orthonormal basis $e_1,...,e_n$
of $\mathcal V$ a  tensor $K$ has the form
\begin{equation}
\begin{array}{rcl}
&&K(e_1,e_1)=\lambda e_1, \ \ \ K(e_1,e_j)=\mu e_j\\
&& K(e_j,e_j)=\mu e_1, \ \ \ K(e_j,e_i)=0
\end{array}
\end{equation}
for $i\ne j$, $i,j>1$, or, equivalently
\begin{equation}
K(X,Y)= (\lambda -3\mu) g(X,e_1)g(Y,e_1)e_1 +\mu g(X,Y)e_1 +
g(X,e_1)Y+\mu g(Y,e_1)X
\end{equation}
for any vectors $X,Y\in \mathcal V$.
 In particular,  the case where $\lambda =3\mu$ appears on the
 Whitney sphere. For Lagrangian pseudospheres one has $\lambda =2\mu$
  (as in Example \ref{example_k=lambda_2over4}),
  for
 Lagrangian-umbilical submanifolds $\lambda =\mu$ (cf. \cite{Ch}).
Observe that $E=\tr _g K$ is equal to $(\lambda+(n-1)\mu)e_1$ and
consequently $e_1=E/{\Vert E\Vert}$.  Since $K_{e_1}$ restricted to
the orthogonal complement $\mathcal D$ to $e_1$ is  a  multiple of
the identity, the orthonormal  vectors $e_2,..., e_n$ can be chosen
in $\mathcal D$ arbitrary.

 If $X,Y$ are orthonormal vectors in $\mathcal V$ then
 \begin{equation}
k(X\wedge Y)=\mu ^2+\mu(\lambda -2\mu)(x_1^2 +y_1^2),
 \end{equation}
 where $X=x_1e_1 +X'$, $Y=y_1e_1+Y'$ for $X', Y'\in\mathcal D$.
 Observe that $x_1^2 +y_1^2\le 1$. Indeed, we have $1=x_1^2
 +\varepsilon _1$, $1=y_1^2+\varepsilon _2$ and $x_1^2y_1^2
 =g(X',Y')^2\le \varepsilon _1\varepsilon _2$,
 where $\varepsilon _1=\Vert X'\Vert ^2$, $\varepsilon _2=\Vert Y'\Vert
 ^2$.
The last inequality is equivalent to $(1-\varepsilon
_1)(1-\varepsilon _2) \le \varepsilon  _1\varepsilon _2$. Hence
$\varepsilon _1 +\varepsilon _2\ge 1$. We now have
$2=x_1^2+y_1^2+(\varepsilon_1 +\varepsilon _2)\ge x_1^2+y_1^2 +1$,
which implies $x_1^2+y_1^2\le 1$.

One now sees that if $\mu(\lambda -2\mu)\ge 0$ then  $\mu^2\le
k(\pi)\le \mu (\lambda-\mu)$ for any vector plane $\pi$ in $\mathcal
V$. Similarly, if $\mu(\lambda -2\mu)\le 0$ then  $\mu(\lambda-\mu)
\le k(\pi)\le\mu ^2$. In particular, if $\lambda =3\mu$ then $\mu
^2\le k(\pi)\le 2\mu^2$. If $\lambda =2\mu$ then $k(\pi)=\mu^2$ (as
in Example \ref{example_k=lambda_2over4}), if $\lambda =\mu$ then
$0\le k(\pi)\le \mu ^2$. If $\lambda =0$ then $-\mu ^2\le k(\pi)\le
\mu ^2$.

}\end{example}

\bigskip

 Denote by  $S^1$ the unit sphere $\{ X\in \mathcal V;\ \ g(X,X)=1\}$ and  by $\Phi$
the  function
$$\Phi: S^1\ni X\to C(X,X,X)=g(K(X,X),X)\in \R.$$
The function $\Phi$ attains its global maximum on $S^1$. This
maximum is non-negative and equals $0$ if and only if $K=0$ on
$\mathcal V$. But $\Phi$ may attain also local  extrema on $S^1$. A
local maximal value can be  non-positive, see Example
\ref{negative_minimum} below.

For orthonormal $U,W\in S^1$  we consider the mapping
 $\Phi(t)=\Phi( \cos t \, U+\sin t\, W)$. Then $\Phi(0)=\Phi (U)$,
$$\Phi'(0) = 3C(U,U,W),$$ $$\Phi ''(0)=3[2C(W,W,U)-C(U,U,U)]$$ and $$\Phi '''(0)=3[-7C(W,U,U)+2C(W,W,W)].$$
Hence if
 $U\in S^1$ is  a point where $\Phi$ attains its (maybe local) maximum and $W\in S^1$ is orthogonal to $U$ then
 \begin{equation}\label{I+IIpochodna}C(U,U,W)=0, \ \ \ \ \ 2C(W,W,U)-C(U,U,U)\le 0\end{equation}
 and, if the equality holds in the last formula then $\Phi'''(0)=0$ and consequently $C(W,W,W)=0$.

The easiest situation which should be taken into account is when the
sectional $K$-curvature is constant for all vector  planes in
$\mathcal V$. In this respect we have

\begin{lemma}\label{ejiri} Let  $g$ be a scalar product on an $n$-dimensional vector space $\mathcal V$.
Let  $K$ be a symmetric $(1,2)$-tensor on $\mathcal V$
 symmetric relative to $g$. If the sectional $K$-curvature is constant and equal to $A$
on $\mathcal V$
then there is an orthonormal basis $e_1,...,e_n$ of  $\mathcal V$ such that
\begin{equation}
 K(e_1,e_1)= \lambda _1e_1,\ \  K(e_1,e_i)=\mu_1 e_i
\end{equation}
\begin{equation}\label{2wEjiri}
K(e_i,e_i)=\mu _1e_1 +...+\mu _{i-1}e_{i-1}+\lambda _i e_i,
\end{equation}
for $i=2,...n$ and
\begin{equation}\label{3wEjiri}
K(e_i,e_j)= \mu _{i}e_j
\end{equation}
for some numbers $\lambda _i$, $\mu _i$ for  $i=1,..., n-1$ and
$j>i$. Moreover
\begin{equation}\label{4wEjiri}
\mu _i= \frac{\lambda_i -\sqrt{\lambda _i ^2
-4A_{i-1}}}{2},\end{equation}
\begin{equation}\label{5wEjiri}
A_i=A_{i-1}-\mu _{i}^2,
\end{equation}
for $i=1,..., n-1$ where $A_0=A$.

If additionally $\tr _g K=0$ then $A\le 0$, $\lambda _i$  and $\mu _i$ are expressed as follows
\begin{equation}\label{wzory_na_lamba_mu}
\lambda _i=(n-i) \sqrt{\frac{-A_{i-1}}{n-i+1}},\ \ \
\mu_i=-\sqrt{\frac{-A_{i-1}}{n-i+1}}.
\end{equation}
In particular, in the last case the numbers $\lambda _i$, $\mu _i$
depend only on $A$ and the dimension of $\mathcal V$. Moreover, if
$A<0$ then $\lambda _i\ne 0$ and $\mu _i\ne 0$ for  every $i$.
\end{lemma}

\proof Let $e_1\in S^1$ be a point where
    $\Phi$ attains its maximum (maybe local). Then

\begin{equation}\label{I_pochodna}
g(K(e_1,e_1), U)=0
\end{equation}
and
\begin{equation}\label{II_pochodna}
 2g(K(e_1,U),U)-g(K(e_1,e_1),e_1)\le 0.\end{equation}
for each vector $U\in S^1$  orthogonal to $e_1$. The subspace
$\{e_1\}^\bot$ is $K_{e_1}$-invariant, hence there is an orthonormal
basis $e'_1,...,e'_n$ of $\mathcal V$ diagonalizing  $K_{e_1}$ such
that $e'_1=e_1$. Let  $\lambda ' _1=\lambda _1$, $\lambda'_2$,...,
$\lambda ' _n$ be eigenvalues corresponding to the eigenvectors
 $e'_1,..., e'_n$ of
$K_{e_1}$. Taking in (\ref{constant_curvature_2}) $X=Z=e'_1$,
$Y=W=e'_i$ for $2\le i\le n$ we get
\begin{equation}
-A+\lambda _1\lambda '_i -(\lambda '_i)^2=0.
\end{equation}
If we regard this equality as an equation relative to $\lambda '_i$,
we obtain at most two  possible values $ \lambda '_i=\frac{\lambda
_1 \pm\sqrt{ \lambda _1^2 -4A}}{2} $.  By  (\ref{II_pochodna}) we
have $2\lambda '_i\le \lambda _1$. Therefore we may exclude the
value $\frac{\lambda _1 +\sqrt{ \lambda _1^2 -4A}}{2}$.  Set
\begin{equation}\mu _1=\frac{\lambda _1 -\sqrt{ \lambda _1^2
-4A}}{2}.\end{equation}

Note that under condition $\lambda_1\ge 0$, if $A<0$ then $\mu_1<0$,
if $A=0$ then $\mu_1=0$ and if $A>0$ then $\mu_1>0$.

Let $\mathcal D=\{e_1\}^\bot$. Vectors belonging to this subspace will be denoted by $X', Y'$ etc.
Denote by  $K'$  the tensor on $\mathcal D$ defined as
 \begin{equation}\label{K'}
 K'=P \circ
K_{| \mathcal D\times \mathcal D},\end{equation} where $P $ is the
orthogonal projection onto $\mathcal D$. Note that this tensor has
the same properties as  $K$. First, it is symmetric and symmetric
relative to  $g$. Moreover
\begin{eqnarray*}
&&A'(g(X',Z')g(Y',W')-g((X',W')g(Y',Z'))\\
&&\ \ \ \ =g(K(X',Z'),K(Y',W'))-g(K(X',W'),K(Y',Z')),
\end{eqnarray*}
for  $X',Y',Z',W'\in \mathcal D$  and some number $A'$. Indeed, using
(\ref{constant_curvature_2})  and the fact that $K(X',Z')=K' (X',Z')+\mu_1 g(X',Z')e_1$ we obtain

\begin{eqnarray}
&&A(g(X',Z')g(Y',W')-g((X',W')g(Y',Z'))\\
&&\ \ \ \ =g(K(X',Z'),K(Y',W'))-g(K(X',W'),K(Y',Z'))\\
&&\ \ \ \ \  =\mu _1^2(g(X',Z')g(Y',W')-g((X',W')g(Y',Z'))\\
&&\ \ \ \ \ \ \ + g(K'(X',Z'),K'(Y',W'))-g(K'(X',W'),K'(Y',Z')).
\end{eqnarray} Thus
 \begin{equation}\label{A'}
 A'= A-\mu _1^2.
 \end{equation}
 In particular, if
$A$ is negative then so is $A'$. Observe also that if the tensor $K$
is traceless then so is $K'$. Indeed, one has the following
equalities
\begin{eqnarray*}
\sum _{i=2}^nK'(e'_i,e'_i)&=& g(\sum _{i=2}^nK'(e'_i,e'_i),
e'_2)e'_2+...+g(\sum _{i=2}^nK'(e'_i,e'_i), e'_n)e'_n\\
&&=g(\sum _{i=2}^nK(e'_i,e'_i),
e'_2)e'_2+...+g(\sum _{i=2}^nK(e'_i,e'_i), e'_n)e'_n\\
&&= -g(K(e'_1,e'_1),e'_2)e'_2-...-g(K(e'_1,e'_1),e'_n)e'_n\\
&&=-g(\lambda _1 e'_1,e'_2)e'_2-...-g(\lambda _1e'_1,e'_n)e'_n=0.
\end{eqnarray*}

We can now apply the  consideration from the beginning of the proof
to the tensor $K'$ on $\mathcal D$. When replacing the basis
$e'_2,...,e'_n$ by a new basis which is  adapted to $K'$ as in the
first part  of the proof we use the fact that ${K_{e_1}}_{\mid
\mathcal D}$ is proportional to the identity. Using then the
induction we get formulas (\ref{2wEjiri})-(\ref{5wEjiri}).

Assume now that $K$ is traceless.
Observe that in this case $A<0$ (if $K\ne 0$).
Namely,  take $Y=Z$  and the trace relative to $g$ in
(\ref{constant_curvature_2}) at places of  $X$ and  $W$. We get
the equality
\begin{equation}
A(1-n)g(Y,Y)=g(K_Y,K_Y),
\end{equation}
for every $Y$
which shows that $A\le  0$ and $A=0$ if and only if $K=0$.

We have the following equalities characterizing  $\lambda _1$ and
$\mu_1$
$$ \mu _1 =\frac{\lambda _1-\sqrt{\lambda _1^2 -4A}}{2},\ \ \ \
(n-1)\mu_1 +\lambda _1=0.$$ Hence
$$\lambda _1=(n-1)\sqrt{\frac{-A}{n}}, \ \ \ \ \  \mu_1
=-\sqrt{\frac{-A}{n}}.$$ By induction we obtain formulas
(\ref{wzory_na_lamba_mu}).\koniec

\begin{remark}\label{remark_ejiri}
{\rm The expression for $K$ in the above proof is obtained in the
following way. The vector $e_1$ is any vector at which $\Phi$
attains a local maximum on $S ^1$, $e_2$ is any  unit vector at
which $\Phi_{|{\mathcal D}\cap S^1}$ attains its local maximum, etc.
We construct a sequence $\lambda _1$, $\mu _1$, $A_1$, $\lambda _2$,
$\mu _2$, $A_2$ etc. For a given $K$ the expression as in the above
lemma is not unique in general. If $K$ is traceless, however, then
the values  $\lambda _i$ and $\mu _i$ are uniquely given. In
particular, if $\Phi$ attains a local maximum on $S^1$ then its
value $\lambda _1$ is equal to $(n-1)\sqrt{\frac{-A}{n}}$. Hence any
local maximum of $\Phi$ on $S^1$ is its global maximum. The same
deals with $\lambda_i$ for $i=2,..., n$.}

\end{remark}

Using the above proof one also gets

\begin{corollary}\label{[K,K]=0}
 If in the above lemma $A=0$, that is, $[K,K]=0$,
  then there is an orthonormal basis $e_1,...,e_n$ of $\mathcal V$ such that
 \begin{equation}\label{postac_dla_[K,K]=0}
  K(e_i,e_i)=\lambda _i e_i ,\ \ \ K(e_i,e_j)=0
 \end{equation} for $i,j=1,...n$ and $i\ne j$. If $[K,K]=0$ and
 $\tr_g K=0$ then $K=0$.
 If $K$ has expression {\rm (\ref{postac_dla_[K,K]=0} )} then
 $[K,K]=0$.
\end{corollary}

In what follows we use some conventions established in the proof of
Lemma \ref{ejiri}. In particular, if some   unit vector $e_1$ is
fixed then the orthogonal complement to $e_1$ in $\mathcal V$ will
be denoted by $\mathcal D$ and $K'$ will be defined by (\ref{K'}).
Moreover by a maximum we shall mean a local maximum unless otherwise
stated.

\begin{lemma}\label{non-positive_curvature}
Let  the sectional $K$-curvature on $\mathcal V$ be  non-positive
for every plane in $\mathcal V$ and  $e_1\in S^1$ be  a point where
$\Phi$ attains a maximum $\lambda_1\ne 0$ on $S^1$. Then the
sectional $K'$-curvature on $\mathcal D$ is also non-positive.
\newline
If moreover the sectional $K$-curvature is negative on $\mathcal V$
then the sectional $K'$-curvature on $\mathcal D$ is negative and
strictly smaller than the $K$-sectional curvature on $\mathcal D$.
\end{lemma}

\proof We have an orthonormal basis $e_2,...,e_n$ of $\mathcal D$
such that $e_1,e_2,...,e_n$ is an orthonormal basis of eigenvectors
of $K_{e_1}$. Let $\lambda _1,..., \lambda _n$ be the corresponding
eigenvalues of $K_{e_1}$.
 We have  $2\lambda _j\le \lambda _1$ for $j>1$.
 Thus if $\lambda _1<0$ then $\lambda _j<0$. If $\lambda _1>0$
  we get $\lambda _j<\lambda _1$.
 By assumption we have
\begin{eqnarray*}
0\ge k(e_1\wedge e_j)&=&g(K(e_1,e_1), K(e_j, e_j))-g(K_{e_1}e_j, K_{e_1}e_j)\\
&&\ \ \
 \ = g(\lambda _1e_1, K(e_j,e_j)) -\lambda _j^2= \lambda _j(\lambda _1-\lambda _j).
\end{eqnarray*}
Therefore $\lambda_j\le 0$ for every $j\ge 2$.


Assume that $\lambda_2,...,\lambda _r$ are non-zero for some $r>1$ and the next eigenvalues  vanish.
We  can define a (positive definite) scalar product $G$ on the space $span \{e_2,..., e_r\}$:
\begin{equation}
  G(X',Y')=-(x_2y_2\lambda_2+...+x_ry_r\lambda _r),
\end{equation}
where $X'=x_2e_2+...+x_re_r, \  Y'=y_2e_2+...+y_re_r$.

Let $X=x_2e_2+...+x_ne_n,\ Y=y_2e_2+...+y_ne_n$ be any two vectors
of $\mathcal D$ and $X',Y'$ be their orthogonal projections onto the
space $span\{ e_2,...,e_r).$ We have
\begin{eqnarray*}K(X,Y)&&=K'(X,Y)+g(K(X,Y),e_1)e_1\\
&&=K'(X,Y) +g(K_{e_1}(x_2e_2+...+x_ne_n), y_2e_2+...+y_ne_n)e_1\\
&&=K'(X,Y) +(x_2y_2\lambda _2+...+x_ny_n\lambda _n)e_1\\
&&= K'(X,Y)-G(X',Y')e_1.
\end{eqnarray*}
Thus
\begin{eqnarray*}
 &&g(K(X,X),K(Y,Y))-g(K(X,Y),K(X,Y))\\
 &&\ \ \ \ = g(K'(X,X),K'(Y,Y))-g(K'(X,Y),K'(X,Y))\\&&
 \ \ \ \ \ \ \ \ \ \ \ \ \ \ \ \ \ \ \ \ \ \ \ \ \ \ \ \ +G(X',X')G(Y',Y')-G(X',Y')^2.
\end{eqnarray*}
Since $G(X',X')G(Y',Y')-G(X',Y')^2$ is non-negative by the Schwarz lemma, we have that
$$g(K'(X,X),K'(Y,Y))-g(K'(X,Y),K'(X,Y))\le 0$$
if  $g(K(X,X),K(Y,Y))-g(K(X,Y),K(X,Y))\le 0$. The above
consideration   provides a  proof of the lemma also in the case
where all the eigenvalues $\lambda _2,..., \lambda _n$ vanish.

If the  sectional $K$-curvature is  negative then $\lambda _j<0$ for
all $j\ge 2$. Thus $G$ is a scalar product on $\mathcal D$.
Therefore, if $X,Y\in \mathcal D$ are orthonormal then $k'(X\wedge
Y)< k(X\wedge Y)$, where $k'$ is the sectional $K'$-curvature.
\koniec

Analogously as above one gets

\begin{lemma}\label{K'-curvature}
 If the sectional $K$-curvature is non-negative or non-positive on
$\mathcal V$ then the sectional $K'$-curvature $k'$ on $\mathcal D$
is not greater than the sectional $K$-curvature. More precisely, if
$\pi$ is a plane in $\mathcal D$ then $k'(\pi)\le k(\pi)$. If the
sectional $K$-curvature is positive or negative on $\mathcal V$ then
$k'(\pi) <k(\pi)$ for every plane in $\mathcal D$.
\end{lemma}
 From the proof of Lemma \ref{non-positive_curvature} we have the
 following useful observation
\begin{lemma}\label{drobiazg}
Let $e_1,...., e_n$ be an orthonormal basis diagonalizing $K_{e_1}$
with corresponding eigenvalues $\lambda_1,...,\lambda _n$. Then
\begin{equation}\label{k(e_1,e_j)}
 k(e_1\wedge e_j)= \lambda _j(\lambda _1-\lambda _j)
\end{equation}
for $j=2,...,n$. If $\lambda _1$ is a maximal value of $\Phi$ on
$S^1$ and $\lambda _1\ge 0$ then $\lambda _j\le \lambda _1$.
 If $\lambda _1>0$ then $\lambda _j<\lambda _1$.
\end{lemma}
\begin{proposition}
Let $\lambda _1$ be a maximal value of $\Phi$ on $S^1$ attained at
$e_1$ and $e_1,..., e_n$ be an eigenbasis of $K_{e_1}$ with
corresponding eigenvalues $\lambda _1,...,\lambda_n$.
 If $\lambda _1=0$ then $k(e_1\wedge e_j)\le 0$ for every $j=2,...,n$. In particular, if
 the sectional $K$-curvature on $\mathcal V$ is positive for all
 planes then $\lambda _1\ne 0$.  If  the structure $(g,K)$
 is trace-free then  $\lambda _1\ge 0$ and $\lambda _1=0$ if and
 only if $K=0$. For a trace-free structure the sectional
 $K$-curvature cannot be non-negative on $\mathcal V$.
\end{proposition}
\proof Assume that $\tr _gK=0$. Then \begin{equation}\label{trace}
\sum_{i=1}^n\lambda _i=0.\end{equation}
 If $\lambda
_1<0$  then, because $2\lambda _j\le \lambda _1$ for $j=2,...,n$, we
have $\lambda _i<0$ for all $i=1,...,n$. This contradicts
(\ref{trace}). If $\lambda _1=0$ then $\lambda _j\le 0$ for $j\ge
2$. By (\ref{trace}) all $\lambda _j=0$. By the remark made after
(\ref{I+IIpochodna}) we have that $K'=0$ and consequently $K=0$.
Suppose that the sectional $K$-curvature is  non-negative on
$\mathcal V$ and $K\ne 0$. Then $\lambda _1>0$. By Lemma
\ref{drobiazg} we have $\lambda _1-\lambda _j>0$. By (\ref{trace})
there is $j>1$ such that $\lambda _j<0$. Hence, using
(\ref{k(e_1,e_j)}),  one gets the contradiction $k(e_1\wedge
e_j)<0$. \koniec

\begin{example}\label{negative_minimum} {\rm It is possible that the sectional $K$-curvature is
positive and $\lambda _1<0$. For instance, define $K$ on the
standard Euclidean space $\R ^2$ with the canonical basis $e_1,e_2$
by
\begin{eqnarray*}
K(e_1,e_1)= -3 e_1,\ \ \  K(e_1,e_2)=-2 e_2,\ \  \ K(e_2,e_2)=
-2e_1.
\end{eqnarray*}
One easily checks (using consideration before Lemma \ref{ejiri})
that $\Phi$ attains a local maximum at $e_1$ and the $K$-curvature
equals $2$.}
\end{example}
 We shall need
\begin{lemma}\label{krzywizny_z_e_1} Let
 $\Phi$ attain its  maximum  $\lambda_1$ at $e_1\in S^1$ and $e_1,...,e_n$ be an
 orthonormal eigenbasis of $K_{e_1}$ with corresponding  eigenvalues $\lambda _i$, $i=1,...,n$.
 If $k(e_1\wedge e_j)<\lambda_1^2/4$ for some $j=2,...,n$ then $2\lambda_j-\lambda_1<0$.
 In particular, if $\lambda_1\ne 0$ and the sectional $K$-curvature is non-positive for all
 planes in $\mathcal V$ then $2\lambda_j-\lambda_1<0$ for every $j=2,...,n$.
\end{lemma}
\proof We know that $2\lambda_j-\lambda_1\le0$ for every
$j=2,...,n$. If
 $2\lambda _j-\lambda_1=0$ then $k(e_1\wedge e_j)=\frac{\lambda _1^2}{4}$.

\koniec

\begin{lemma}\label{k(e_1,X)}
 Let $\lambda_1$ be a  maximal value of $\Phi$ on $S^1$ attained at $e_1\in S^1$. Let $X\in S^1$ be orthogonal
 to $e_1$. Then $k(e_1\wedge X)\le \frac{\lambda_1^2}{4}$ and the equality holds if and only if
 $X$ is an eigenvector of $K_{e_1}$ with eigenvalue $\frac{\lambda_1}{2}$.
\end{lemma}
\proof Assume first that $X\in S^1$ is an eigenvector of $K_{e_1}$
with corresponding eigenvalue $\mu$. Then $k(e_1\wedge X)=-\mu^2
+\mu\lambda _1$. Since the function $\R\ni t\to -t^2+\lambda_1 t$
attains its maximum $\lambda _1^2/4$  for $t=\lambda_ 1/2$, we have
that $k(e_1\wedge X)\le \lambda_1^2/4$ and the equality holds   if
and only if $\mu=\lambda _1/2$.

As usual, let $e_1,...,e_n$ be an orthonormal eigenbasis for
$K_{e_1}$ and $\lambda_1,..., \lambda _n$ be the corresponding
eigenvalues. Let $X=x_2e_2+...+x_ne_n\in S^1$ be orthogonal to $e_1$
but not necessary an eigenvector of $K_{e_1}$. One now gets
\begin{equation}
\begin{array}{lcr}
 &&k(e_1\wedge X)=g(K(e_1,e_1), K(X,X))-g(K(e_1,X), K(e_1,X))\\
 &&\ \ \ \ \ \ \ =\lambda_1g(K_{e_1}(x_2e_2+...+x_ne_n), x_2e_2+...+x_ne_n)\\
 &&\ \ \ \ \ \ \ \ \ \ \ \  \ -g(K_{e_1}(x_2e_2+...+x_ne_n),K_{e_1}(x_2e_2+...+x_ne_n))\\
 &&\ \ \ \ \ \  \ \ \ \ \ \ \ \ \ \
 =\lambda_1(x_2^2\lambda_2+...+x_n^2\lambda_n)-(x_2^2\lambda_2^2+...+x_{n}^2\lambda_n^2)\\
 &&\ \ \ \ \ \ \ \ \ \ \ \ \ \ \ \ \ \ =k(e_1\wedge e_2)x^2_2+...+k(e_1\wedge e_n)x^2_n\\
 &&\ \ \ \ \ \ \ \ \ \ \ \ \ \ \ \ \ \ \ \
 \ \ \ \ \ \ \le \frac{\lambda _1^2}{4}x_2^2+...+\frac{\lambda_1^2}{4}x_n^2=
 \frac{\lambda_1 ^2}{4}.
 \end{array}
\end{equation}
In this formula the equality holds if and only if for each $j=2,...,n$ either $x_j=0$ or $\lambda_j=\lambda_1/2$.
Assume that
$x_2,...,x_r$ are not zero and the next coordinates of $X$   vanish.
Then $\lambda_2=...=\lambda _r=\lambda_1/2$ and one sees that $K_{e_1}X=\frac{\lambda_1}{2}X$.\koniec

\begin{lemma}
    Let $\Phi$ attain its  maximum at $e_1$ and
  $e_1,...,e_n$ be an eigenbasis with corresponding eigenvalues $\lambda _1,...,\lambda _n$.
  If $2\lambda _j-\lambda _1<0$ then for each $X\in S^1$ orthogonal to $e_1$ we have
  $$2C(X,X,e_1)-\lambda _1<0.$$
\end{lemma}
\proof Let $X=x_2e_2+...+x_ne_n$. Then
\begin{eqnarray*}
&&2C(X,X,e_1)=2g(K_{e_1}(x_2e_2+...+x_ne_n), x_2e_2+...+x_ne_n)\\
&&\ \ \ \ \ \  \ =2(\lambda_2x_2^2+...+\lambda
_nx_n^2)<\lambda_1x_2^2 +  ... + \lambda _1 x_n^2=\lambda _1.
\end{eqnarray*}
\koniec

From the above lemmas we immediately get
\begin{proposition}
  Let $\lambda_1$ be a  maximal value of $\Phi$ on $S^1$ attained at $e_1\in S^1$ and  $X\in S^1$ be orthogonal
 to $e_1$. Then $k(e_1\wedge X)<\frac{\lambda _1^2}{4}$ if and only if $2C(X,X,e_1)-\lambda _1<0$.
\end{proposition}
By Lemma \ref{k(e_1,X)} we know  that if the $K$-sectional curvature
is constant then its value is less than or equal to $\frac{\lambda
_1^2}{4}$ for any maximal value $\lambda _1$ on $S^1$.

\begin{proposition}
 Assume that the sectional $K$-curvature on $\mathcal V$ is constant and equal to $A=\frac{\lambda_1^2}{4}$
where $\lambda _1$ is a maximal value of $\Phi $ on $S^1$. Then
there is an orthonormal  basis $e_1,...,e_n$ of $\mathcal V$
relative to which $K$
 has expression as in {\rm Example \ref{example_k=lambda_2over4}}.
 \end{proposition}
\proof
 By Lemmas \ref{ejiri} and
  \ref{k(e_1,X)} we have $K_{e_1}X=\frac{\lambda _1}{2}X$ for any $X$ orthogonal to $e_1$.
   Since $\Phi$ attains a maximum $\lambda _1$ at $e_1$ and $2C(X,X,e_1)-\lambda _1=0$, by
  the observation made in the sentence
  containing
 (\ref{I+IIpochodna}) we know that $C(X,X,X)=0$ for every $X$ orthogonal to $e_1$. \koniec

Consider now the vector $E=\tr _gK$. If the sectional $K$-curvature
is constant and equal to $\frac{\lambda_1^2}{4}$, where $\lambda_1$
is a maximum of $\Phi$ on $S^1$, then, by (\ref{trK}),
$E=\frac{n+1}{2}\lambda _1 e_1$. Therefore, if we have a statistical
structure $(g,K)$ on  a manifold $M$ of constant sectional
$K$-curvature equal to $\frac{\lambda_1^2}{4}$ at each point of $M$
then $\lambda _1$ is constant and $e_1$ is a smooth vector field on
$M$.

Note that the assumption that $\lambda $ in  Example
\ref{example_k=lambda_2over4} is a maximal value of $\Phi$ is not
needed. We have the following characterizations of the structure
from Example \ref{example_k=lambda_2over4}

\begin{thm}\label{characterization}
Structures in {\rm Example \ref{example_k=lambda_2over4}} are
characterized by the conjunction of the following conditions:
\newline
1) $E$ is an eigenvector of $K_E$,
\newline
2) $K_E$ restricted to the orthogonal complement to $E$ is a
multiple of the identity,
\newline
3) the sectional $K$-curvature on $\mathcal V$ is  a positive
constant $A$,
\newline
4) $\Vert E\Vert =(n+1)\sqrt{A}$.
\end{thm}

\proof Of course, if the structure is as in Example
\ref{example_k=lambda_2over4} then all conditions $1)-4)$ are
satisfied. Assume that the conditions $1)-4)$ are fulfilled. By $1)$
and $2)$ we know that there exist numbers $\lambda $ and $\mu$ such
that $K_{e_1}e_1=\lambda e_1$ and $K_{e_1}e_i=\mu e_i$, for
$i=2,...,n$, where $e_1=\pm E/\Vert E\Vert$ and $e_1,...,e_n$ is an
orthonormal basis of $\mathcal V$. We choose the sign of $e_1$ in
such a way that $\lambda\ge 0$. By 3), similarly as in the proof of
Lemma \ref{ejiri}, we obtain
\begin{equation}\label{lambda_mu}
\mu=\frac{\lambda \pm \sqrt{\lambda^2 -4B^2}}{2},
\end{equation}
where $A=B^2$, for some $B>0$. In particular, we have $\lambda
^2-4B^2\ge 0$, which implies that  $\lambda -2B\ge0$. By $4)$ we
have
\begin{equation}\label{4}
(n-1)\mu +\lambda =(n+1)B.\end{equation} Inserting (\ref{lambda_mu})
into (\ref{4}) we get
$$
\pm(n-1)\sqrt{\lambda -2B}\sqrt{\lambda
+2B}=-(n+1)\sqrt{\lambda-2B}\sqrt{\lambda -2B}.
$$
Assume that $\lambda\ne 2B$. Then $\mp (n-1)\sqrt{\lambda
+2B}=(n+1)\sqrt{\lambda-2B}$ and consequently $ (n-1)^2(\lambda
+2B)=(n+1)^2(\lambda -2B)$. It follows that
$$\lambda =\frac{n^2+1}{n}B.$$
Inserting this into (\ref{lambda_mu}) one gets $\mu=nB$ or $\mu
=B/n$. Using now (\ref{4}) we obtain contradictions. Therefore
$\lambda=2B$ and, by (\ref{lambda_mu}) $\mu =\lambda /2$. It follows
that $A=\lambda ^4/4$. We can now go back to the  proof of Lemma
\ref{ejiri}. By (\ref{A'}) we see that the sectional $K'$-curvature
on $\mathcal D$ vanishes. Hence $K'$ has expression as  in Corollary
(\ref{[K,K]=0}). But  $E$ is proportional to $e_1$, hence $K'=0$ and
consequently $K$ has expression as in Example
\ref{example_k=lambda_2over4}.\koniec

\begin{thm} Let $(g,K)$ be a statistical structure on $M$
such that at each point $p$ of $M$ the tensor $K_p$ is as in Theorem
{\rm \ref{characterization}}.
If $\hat\nabla K$ is symmetric and $\div E$ is constant then  the
sectional curvature (for $g$) by any plane containing $E$ is
non-positive.
 If
$\hat\nabla E=0$ then  $\hat\nabla K=0$ on $M$.
\end{thm}
\proof We can assume that $M$ is connected. Since $\hat \nabla K$ is
symmetric, the sectional $K$-curvature is constant on $M$. We have
$E=\Lambda e_1$ where $\Lambda $ is a smooth function and $e_1$ is a
smooth unit vector field on $M$.  $\lambda$ is a constant function
on $M$. Locally we can extend $e_1$ to a smooth orthonormal frame
$e_1,...,e_n$. In such a frame $K$ has expression as in Example
\ref{example_k=lambda_2over4}. Then $\Lambda =(n+1)\lambda/2$. Let
$\hat \nabla _{e_i}e_j=\sum _{k=1}^n\omega^k_j(e_i)e_k$.

By a  straightforward computation one gets for  mutually different
$i,j,l\ge 2$
\begin{equation}\label{nablaK}
\begin{array}{rcl}
&&(\hat\nabla _{e_1}K)(e_1,e_1)=0\\
&&(\hat\nabla _{e_i}K)(e_i,e_i)=\frac{\lambda}{2}\sum_{k\ne 1}\omega
_1^k(e_i)e_k+\lambda \omega _1^i(e_i)e_i \\
&&(\hat\nabla _{e_i}K)(e_1,e_1)=0\\
&&(\hat\nabla _{e_1}K)(e_i,e_1)=0\\
&&(\hat\nabla _{e_1}K)(e_i,e_i)=\frac{\lambda}{2}\sum_{k\ne
1,i}\omega
^k_1(e_1)e_k +\frac{3}{2}\lambda\omega^i_1(e_1)e_i\\
&&(\hat\nabla _{e_i}K)(e_1,e_i)=0\\
&&(\hat\nabla
_{e_1}K)(e_i,e_j)=\frac{\lambda}{2}\omega^i_1(e_1)e_j+\frac{\lambda}{2}\omega^j_1(e_1)e_i\\
&&(\hat\nabla _{e_i}K)(e_1,e_j)=0\\
&&(\hat\nabla _{e_i}K)(e_j,e_j)=\frac{\lambda }{2}\sum_{k\ne
1}\omega
^k_1(e_i)e_k+\lambda\omega^j_1(e_i)e_j\\
&&(\hat\nabla _{e_j}K)(e_i,e_j)=\frac{\lambda}{2}\omega _1^j(e_j)e_i
+\frac{\lambda}{2}\omega _1^i(e_j)e_j\\
&&(\hat\nabla
_{e_i}K)(e_j,e_l)=\frac{\lambda}{2}\omega_1^j(e_i)e_l+\frac{\lambda}{2}\omega_1^l(e_i)e_j\\
&&(\hat\nabla _{e_j}K)(e_i,e_l)=\frac{\lambda}{2}\omega_1^i(e_j)
e_l+\frac{\lambda}{2}\omega_1^l(e_j) e_i.
\end{array}
\end{equation}
One now sees that if $\hat\nabla K$ is symmetric then $\hat\nabla
_{e_1}e_1=0$,  and $\hat\nabla _{e_i}e_1=\alpha e_i$ for some
function $\alpha$ for every $i=2,...,n$. It implies that
 $$g(\hat
R(e_i,e_1)e_1,e_i)= -(e_1\alpha +2\alpha ^2)$$ for every
$i=2,...,n$.
 Since $\div E
=\frac{n^2-1}{2}\lambda \alpha $\ is constant,  the function
$\alpha$ is constant if $\div E$ is constant. Consequently $g(\hat
R(e_i,e_1)e_1,e_i)=-2\alpha ^2$. If $\hat\nabla E=0$ then
$\hat\nabla e_1=0$ and formulas (\ref{nablaK}) imply $\hat\nabla
K=0$.\koniec

If $J$ is an endomorphism of $\mathcal V$ and $T$ is a tensor on
$\mathcal V$ then $J\cdot T$ will mean that $J$ acts as a
differentiation on $T$. If $\mathcal R$ is a tensor of type $(1,3)$
and $\mathcal R(X,Y)$  denotes the endomorphism determined by
$\mathcal R$ then the equality $\mathcal R\cdot T=0$ means that
$\mathcal R(X,Y)\cdot T=0$ for every $X,Y\in \mathcal V$. If $X\in
\mathcal V$ then $\mathcal RX=0$ means that $\mathcal R(Y,Z)X=0$ for
every $Y,Z\in \mathcal V$. The same convention will be used  for
tensor fields on manifolds.

\begin{lemma}\label{JK=0_1}
 Let $J$ be an endomorphism of $\mathcal V$ such that $J\cdot g=0$, where $J$ is regarded as a differentiation.
 If the sectional $K$-curvature is negative for every plane of $\mathcal V$ and $J\cdot K=0$ then
 $J=0$.
\end{lemma}
\proof  As usual take $e_1\in S^1$ where $\Phi$ attains its maximum
and an orthonormal eigenbasis $e_1,..., e_n$ of $K_{e_1}$ with
corresponding eigenvalues $\lambda _1,...,\lambda _n$. By Lemma
\ref{krzywizny_z_e_1} we know that $2\lambda _i - \lambda _1< 0$ for
all $i=2,..,n$. Using the fact that $J$ is skew-symmetric relative
to $g$ we obtain
\begin{eqnarray*}
 0=(J\cdot K)(e_1,e_1) &=&J(K(e_1,e_1))-2K(Je_1,e_1)\\
  &=& \lambda_1 \sum_{j=2}^n g(Je_1, e_j)e_j -2K (\sum_{j=2}^ng(Je_1, e_j)e_j, e_1)\\
  &=&\sum_{j=2}^n (\lambda _1-2\lambda _j )g(Je_1,e_j)e_j.
\end{eqnarray*}
Using also the fact that $g(Je_1,e_1)=0$, we get $Je_1=0$. In
particular, the orthogonal complement $\mathcal D$ to $e_1$ in
$\mathcal V$ is $J$-invariant. Let $K'$ be given by (\ref{K'}) and
$J'$ stands for the restriction of $J$ to $\mathcal D$.  For
$X',Y'\in \mathcal D$ we get (using the skew-symmetry of $J$, the
condition $J\cdot K=0$ and the equality $Je_1=0$)
\begin{eqnarray*}
&&(J'\cdot K')(X',Y')=J(K(X',Y')-g(K(X',Y'),e_1)e_1)\\
&&\ \ \ \ \ \ \ \ \  -K(JX',Y')+g(K(JX',Y'),e_1)e_1-K(X',JY')+g(K(X',JY'),e_1)e_1\\
&&\ \ \ \ \ \ \ \ \ \ \ \ \ \ \  \ \ \ \ \ \  \ =(J\cdot
K)(X',Y')+g(J(K(X',Y')),e_1)e_1=0.
\end{eqnarray*}
 By Lemma \ref{K'-curvature} we see that $K'\ne 0$  and the sectional $K'$-curvature
 is negative on $\mathcal D$. We can now apply  the same as above arguments for the objects $K', J'$ on
  $\mathcal D$ and continue the proof using induction.\koniec

Using Lemma \ref{krzywizny_z_e_1} and the first part of the proof of
Lemma \ref{JK=0_1} we obtain

\begin{lemma}
 Let $J$ be an endomorphism of $\mathcal V$ such that $J\cdot g=0$, where $J$ is treated as a differentiation.
 Assume that $\lambda  _1\ne 0$ is a maximal value of $\Phi$ on $S^1$ attained at $e_1\in S^1$.
 If the sectional $K$-curvature is smaller  than $\lambda_1^2/4$
 for every plane of $\mathcal V$ and
 $J\cdot K=0$
 then $Je_1=0$.
 \end{lemma}

\begin{thm}\label{non-positive_curvature}
 If the sectional $K$-curvature is non-positive on $\mathcal V$ and
  $[K,K]\cdot K=0$ then the sectional $K$-curvature vanishes on $\mathcal V$.
 \end{thm}
\proof We  can modify the proof of Lemma \ref{JK=0_1}. Assume that
$K\ne 0$. Let $\lambda _1>0$ be  a maximal value of $\Phi$ on $S^1$
attained at $e_1$ and $e_1,...,e_n$ be an eigenbasis of $K_{e_1}$
with corresponding eigenvalues $\lambda _1,..., \lambda _n$. By
Lemma \ref{krzywizny_z_e_1} we have $2\lambda _j-\lambda _1\ne 0$.
As in the proof of Lemma \ref{JK=0_1}  we get $[K,K]e_1=0$. It
follows that $\lambda_ j(\lambda _1-\lambda_j)=0$ for every $j\ge
2$, hence $\lambda_j=0$ for  every $j\ge2$ (because $\lambda
_j<\lambda _1$ if $\lambda _1>0$, see Lemma \ref{drobiazg}). Thus
$K(e_1,X')=0$ and consequently $g(e_1, K(X',Y'))=0$ for every $X',
Y'\in \mathcal D$.
Let $X,Y$ be any vectors of $\mathcal V$ and $X=x_1 e_1 +X', Y=y_1
e_1 +Y'$ for some $X',Y'\in \mathcal D$. We now have
\begin{eqnarray*}
&&g(K(X,X),K(Y,Y))-g(K(X,Y), K(X,Y))\\
&&=g(x_1 ^2\lambda _1e_1+K(X',X'),y_1^2\lambda_1 e_1+K(Y',Y'))\\&&\
\ \ \ \ \
-g(x_1y_1\lambda_1e_1+K(X',Y'),x_1y_1\lambda_1e_1+K(X',Y'))\\
&&=g(K'(X',X') +g(K(X',X'), e_1)e_1, K'(Y',Y')+g(K(Y',Y')e_1)e_1)\\
&&\ \ \ \ \ \
\ -g(K'(X',Y')+g(K(X',Y')e_1)e_1,K'(X',Y')+g(K(X',Y')e_1)e_1)\\
&&=g(K'(X',X'),K'(Y',Y'))-g(K'(X',Y'), K'(X',Y')),
\end{eqnarray*}
that is,
\begin{equation}
\begin{array}{rcl}
&&g(K(X,X),K(Y,Y))-g(K(X,Y), K(X,Y))\\&&\ \ \
 \ =g(K'(X',X'),K'(Y',Y'))-g(K'(X',Y'), K'(X',Y')).
 \end{array}
\end{equation}
  It follows that the sectional $K$-curvature vanishes on $\mathcal
V$ if $K'=0$ and if   $K'\ne 0$ the sectional  $K'$-curvature on
$\mathcal D$ is non-positive.  In the last case we argue  as above
for the structure $K'$ on $\mathcal D$. We get $[K',K']e_2=0$ and
$K'(e_2,X'')=0$ for every $X''\in \mathcal D$ orthogonal to $e_2$,
where $e_2\in S^1\cap \mathcal D$ is a point where $\Phi_{|S^1\cap
\mathcal D} $ attains its positive maximal value. Then we continue
the proof by induction using the same type of arguments as above and
we obtain the expression for $K$ as in Corollary \ref{[K,K]=0}.
\koniec

As  consequences of  Theorem \ref{non-positive_curvature} we obtain
\begin{corollary}
If $(g,K)$ is a Hessian structure on $M$ with non-negative sectional
curvature of $g$ and such that $\hat R\cdot K=0$ then  $g$ is
flat.\end{corollary} Lemma \ref{JK=0_1} yields
\begin{thm}
 If $(g,K)$ is a statistical structure on a manifold  $M$, the sectional $K$-curvature  is negative on $M$
 and  $\hat R\cdot K=0$ then  $g$ is flat.
\end{thm}
In the following theorem $M$ is $n$-dimensional and the complex
space form has complex dimension $n$.
\begin{thm}
 If $M$ is a totally real submanifold of  the complex space form of holomorphic sectional curvature $4c$,
 the sectional curvature  of $M$ is smaller than $c$ on $M$ and  $\hat R\cdot K=0$, where  $K$ is
 the second fundamental tensor of the submanifold then  $\hat R=0$.
\end{thm}
\proof We have the  following Gauss equation for a totally real
submanifold in the complex space form
\begin{equation}
c[g(Y,Z)X-g(X,Z)Y]=\hat R(X,Y)Z-[K_X,K_Y]Z
\end{equation}
for every $X,Y,Z$ tangent to $M$. Hence the sectional $K$-curvature
equals to the difference between the sectional curvature for $g$ and
$c$.  Therefore the assumption of the theorem  says that the
sectional $K$-curvature is negative so we can use Lemma
\ref{JK=0_1}.\koniec

\begin{thm}
 Let $(g,K)$ be a statistical structure on  a connected manifold $M$  and $g$ has constant sectional curvature.  If
 at some point $p$ of $ M$ the equality $\hat R \cdot K=0$ holds and  the sectional $K$-curvature is positive
 on $T_pM$
 and  strictly smaller
 than the maximal value of $\Phi$  on the unit sphere in $T_pM$ then $g$ is a flat metric.
 \end{thm}
\proof By the proof of Lemma \ref{JK=0_1} we have that $\hat
R_pe_1=0$. Hence $\hat R=0$.\koniec

\section{A smoothness lemma and its consequences for statistical structures of constant sectional
$K$-curvature}
 Although for a statistical structure $(g,K)$ with
constant $K$-sectional curvature on a manifold $M$ at each point of
$M$ we can find an orthonormal frame  for which $K$ has  expression
as in Lemma \ref{ejiri}, it is not possible, in general, to find a
smooth orthonormal local frame relative to which $K$ has this nice
expression. Even to find  a smooth local vector field $e_1$ at which
$\Phi$ attains a maximum  makes a problem. We shall now prove (see
Lemma \ref{smoothness_lemma} below) that in some cases it is
possible. Since in the proof we shall use the multiple Lagrange
method, we can only get a vector field at which $\Phi$ attains a
local maximum (even if we start with a global maximum at some point
$p\in M$). It is why we have used local maxima in our
considerations, for instance in Lemma \ref{ejiri}. Since the author
of this paper was unable to find  references for Lemma
\ref{smoothness_lemma} with a rigorous proof, we provide a detailed
proof. We shall use  Lemma \ref{smoothness_lemma}  only for the
cubic form of  statistical structures, but we formulate and prove
the result for symmetric forms of any degree.

We shall start  with a topological lemma
\begin{lemma}\label{topological_lemma}
 Let $\pi :\mathcal H\to M$ be a locally trivial bundle with a compact standard fiber $ H$
and let $\psi:\mathcal H\to  T$ be a continues mapping into a
topological space
 $ T$. If $\mathcal H_p\subset
\psi ^{-1}(B)$ for some open subset $B\subset T$  then there is a
neighborhood $\mathcal  U$ of $p$ in $M$ such that $\bigcup _{x\in
\mathcal U}\mathcal H_x\subset \psi ^{-1}(B)$.\end{lemma}

\proof  We can assume that in some neighborhood $M'$ of $p$ the
bundle of the shape $ M' \times H$. For every $v\in \mathcal
H_p=\{p\}\times H$ there is a neighborhood $\mathcal U_v$ of $v$ in
$\mathcal H$  such that $ \psi (\mathcal U_{v}) \subset B$. We can
assume that $\mathcal U _v= U_v\times D_{v}$, where $D_{v}$ is an
open subset in  $H$ and $U_v$ is an open neighborhood of $p$. Of
course $\bigcup _{v\in \mathcal H_p}D_{v}$ contains $H$. We choose a
finite subcovering  $ D_{v_1},..., D_{v_r}$  of  $H$ and take $
\mathcal U=\bigcap _{i=1}^{r} U_{v_i}$. Let $(x,v)\in \mathcal U
\times H$. Then $(x,v)\in\mathcal U_{v_i}$ for some $i=1,...,r$.
Hence $\psi (U)\subset B$. \koniec

By a Riemannian vector bundle we mean a vector bundle $\mathcal W\to
M$ for which each fibre $\mathcal W_p$ has a scalar product $g_p$
and the assignment $ p\to g_p$ is smooth.
\begin{lemma}\label{smoothness_lemma}
 Let $\mathcal W$ be a Riemannian vector bundle over $M$ and $U\mathcal W$ be its unit sphere bundle.
  Assume that $C$ is a smooth field of symmetric
 $(0,k)$-tensors on  $\mathcal W$ and $\Phi$ is
 defined as follows $\Phi: U\mathcal W\ni X\to C(X,...,X)\in \R$. Assume that at each point $p\in M$
 the function $\Phi_p=\Phi_{|U\mathcal W _p}$ has the following property:
 \newline
 {\rm (*)}
 If $\Phi _p$ attains its local maximum on $U\mathcal W_p$ at $X_0$ then
 $(k-1)C(U,U,X_0,...,X_0)-C(X_0,...,X_0)\ne 0$ for every $U \in U\mathcal W _p$ orthogonal to $X_0$.
 \newline
 Then for every $p\in M$ there is a smooth unit section  $e_1$  of $\mathcal
 W$,
 defined in some neighborhood of $p$, such that
 $\Phi _x$ attains its (local) maximum on $U\mathcal W_x$  at
$e_1(x)$  for each  $x$ from this neighborhood.
\end{lemma}

\proof Let   $p\in M$  be a fixed point. Denote by $n$ the rank of
the bundle $\mathcal W$. Let $e_1$ be a point of $U\mathcal W_p$ at
which $\Phi$ attains a local maximum. Then
\begin{equation}
 C(U,e_1,...,e_1)=0
\end{equation}
and
\begin{equation}
 (k-1)C(U,U,e_1,...,e_1)-C(e_1,...,e_1)<0
\end{equation}
for any $U\in U\mathcal W_p$ orthogonal to $e_1$. Let $G$ be a
symmetric 2-form on $\mathcal W_p$ given by $G(X,Y)=
C(X,Y,e_1,...,e_1)$. Since $G(U,e_1)=0$ for every $U$ orthogonal to
$e_1$ there is an orthogonal basis $e_1, ..., e_n$ of $\mathcal W_p$
diagonalizing $G$. Let  $\lambda _1$,..., $\lambda _n$ be
eigenvalues corresponding to the  basis $e_1,...,e_n$. We have
$(k-1)\lambda_i-\lambda _1<0$ for $i=2,...,n$. Extend the
orthonormal frame to any local orthonormal frame, say $E_1,...,E_n$
in a neighborhood $\mathcal U$ of $p$. Let $C_{i_1...i_k}$ be the
coordinates of the  form $C$  relative to this local frame, that is,
$C_{i_1...i_k}=C(E_{i_1},..., E_{i_k})$. Consider functions
$f:\mathcal U\times \R^n\times \R\to \R$ defined as follows
$$f(x,y_1,...,y_n, \lambda)=\sum _{i_1,...,i_k=1}^nC_{i_1...i_k}(x)y_{i_1}\cdot\cdot\cdot y_{i_k}
-\lambda (y_1^2+...+y_n^2-1).$$ By the Lagrange method one knows
that at a fixed point $x\in \mathcal U$, the extrema of the function
$\sum _{i_1,...,i_k=1}^nC_{i_1...i_k}(x)y_{i_1}\cdot\cdot\cdot
y_{i_k}$  on the sphere $y_1^2+...+y_n^2-1=0$ are in the set
described by the system of equations

\begin{eqnarray*}
\frac{\partial f}{\partial y_1}&=&k\sum_{i_2,...,i_k=1}^n
C_{1i_2...i_k}(x)y_{i_2}\cdot\cdot\cdot y_{i_k}-2\lambda
y_1=0\\
\frac{\partial f}{\partial y_2}&=&k\sum_{i_2,..., i_k=1}^n
C_{2i_2...i_k}(x)y_{i_2}\cdot\cdot\cdot y_{i_k} -2\lambda
y_2=0\\
&&.\\
&&.\\
&&.\\
\frac{\partial f}{\partial y_n}&=&k\sum_{i_2,...,i_k=1}^n
C_{ni_2...i_k}(x)y_{i_2}\cdot\cdot\cdot y_{i_k} -2\lambda
y_n=0\\
\frac{\partial f}{\partial \lambda}&=&-(y_1^2+...+y_n^2-1)=0
\end{eqnarray*}
Define the functions
\begin{equation}
F_i(x,y_1,...y_n,\lambda)=k\sum_{i_2,..., i_k=1}^n
C_{ii_2...i_k}(x)y_{ i_2}\cdot\cdot\cdot y_{i_k}-2\lambda y_i
\end{equation}
for  $i=1,...,n$ and $$F_{n+1}(x,y_1,...,y_n,\lambda
)=y_1^2+...+y_n^2-1.$$ Set $y_{n+1}=\lambda$. Let $F=(F_1,...,
F_{n+1}): \mathcal U\times \R ^{n+1}\to \R^{n+1}$.

We want to find smooth functions $y_1(x),..., y_n(x), \lambda (x)$,
which satisfy the equation $F(x, y_1(x),..., y_{n+1}(x),\lambda
(x))=(0,...,0)$  and satisfy the initial conditions
$y_1(p)=1,y_2(p)=0..., y_n(p)=0$, $y_{n+1}(p)=\lambda
(p)=\frac{k}{2}C_{1...1}=\frac{k}{2}\lambda _1$. The initial
conditions follow from the fact that the vector  $e_1=(1,0,...,0)$
is among solutions of the above system of equations and
 $\lambda (p)$  can be computed from the first equation of the system.
We shall now use the implicit function theorem. To this aim, it is
sufficient to check that
$$\det \left (\frac{\partial F_i}{\partial y_j}\right )(p, e_1, \frac{k}{2}\lambda_1)\ne 0.$$
We have
$$\frac{\partial F_i}{\partial y_j} =k(k-1)\sum _{i_3,...,i_k=1}^n
C_{iji_3...i_k}(x)y_{i_3}\cdot\cdot\cdot y_{i_k}-2\delta
_{ij}\lambda
$$ for $i,j=1,...,n$. It follows that at the initial values
 we have
 $$\frac{\partial F_i}{\partial y_j}(p,(1,0,...,0),\frac
 {k}{2}\lambda _1)= k[(k-1)C_{ij1...1}(p)-\delta_{ij}\lambda_1]$$
for $i,j=1,...,n$. In particular
$$\frac{\partial F_1}{\partial y_1}(p,(1,0,...,0),\frac
 {k}{2}\lambda _1)=k[(k-1)C_{1...1}(p) -\lambda_1]=k(k-2)\lambda _1,$$
$$\frac{\partial F_i}{\partial y_i}(p,(1,0,...,0),\frac
 {k}{2}\lambda _1)= k[(k-1)C_{ii1...1}(p)-\lambda
_1)]=k[(k-1)\lambda_i-\lambda _1],$$
$$\frac{\partial F_i}{\partial y_j}(p,(1,0,...,0),\frac
 {k}{2}\lambda _1)=0$$
for $i\ne j$, $i,j=2,... ,n$. Moreover
$$\frac{\partial F_1}{\partial y_{n+1}}(p,(1,0,...,0),\frac
 {k}{2}\lambda _1)=-2,$$
 $$\frac{\partial F_j}{\partial y_{n+1}}(p,(1,0,...,0),\frac
 {k}{2}\lambda _1)=0$$
 for $j=2,...,n$;
$$\frac{\partial F_{n+1}}{\partial y_{1}}(p,(1,0,...,0),\frac
 {k}{2}\lambda _1)=2,$$
$$\frac{\partial F_{n+1}}{\partial y_{i}}(p,(1,0,...,0),\frac
 {k}{2}\lambda _1)=0$$
 for $i=2,...,n$;
$$\frac{\partial F_{n+1}}{\partial y_{n+1}}(p,(1,0,...,0),\frac
 {k}{2}\lambda _1)=0.$$
It is now clear that $\det \left (\frac{\partial F_i}{\partial
y_j}\right )(p, e_1,\frac{k}{2}\lambda_1)\ne 0.$

Let $y_1(x),...,y_n(x), \lambda (x)$ be the solution of of our
implicit function problem. Denote by  $e_1$ the section of $\mathcal
W$ given by $y_1E_1+...+y_nE_n$. Since the condition
$F(x,e_1(x),\lambda (x))=0$ is satisfied, at each point of some
neighborhood $\mathcal U'$ of $p$, we have that $C(U,e_1,...,e_1)=0$
for every $U$
 orthogonal to $e_1$, $U\in U\mathcal W_x$ at each $x\in \mathcal U'$. To see this it is sufficient to multiply
 each $F_i(x,e_1(x),\lambda (x))$
 by $U_i$ (where $U=U_1E_1+...+U_nE_n$) and make  summation relative to $i=1,...,n$.

 Using now Lemma \ref{topological_lemma} one sees that since
 $C(e_1,...,e_1)>2C(U,U,e_1,...,e_1)$ for each $U\in U\mathcal W_p$, there is a neighborhood
  $\mathcal U''\subset \mathcal U'$ of $p$ such that
 $C(e_1,...,e_1)>2C(U,U,e_1,...,e_1)$  for each $U\in U\mathcal W_x$  and $x\in
 \mathcal U''$.
 Indeed, it is sufficient to take as $\mathcal H$ the bundle $U\mathcal W_{|{\mathcal U'}}\cap
 \mathcal D$, where $\mathcal D$ is  the orthogonal complement to
 $e_1$ in the  bundle $\mathcal W_{| \mathcal U'}$  and define $\psi$
 as the mapping
 $$\psi: \mathcal H \ni V\to C(V,V,e_1,...,e_1)\in \R.$$

It follows that for every $x\in \mathcal U''$ the mapping  $\Phi _x$
attains at $e_1(x)$  a local maximum.
\koniec

\begin{lemma}\label{nablaKsymetryczne}
Let $(g,K) $ be a statistical structure on a manifold $M$ and its
sectional $K$-curvature is constant. Assume that  for each point $p$
of $M$ there is a local orthonormal frame $e_1,...,e_n$  around $p$
relative to which $K$ has expression as in {\rm Lemma \ref{ejiri}};
$\lambda _i$, $\mu _i$ are constant and $\lambda _i-2\mu _i\ne 0$
for every $i=1,...,n-1$. If $\hat\nabla K$ is symmetric then
$\hat\nabla e_j=0$ for every $j=1,...,n$. In particular, $\hat R=0$
and $\hat\nabla K=0$ on $M$.
\end{lemma}
\proof For every $j>1$ we have
$$g((\hat\nabla _{e_j}K)(e_1,e_1), e_1)=0,\ \ \ \
g((\hat\nabla_{e_1}K)(e_j,e_1), e_1)=(\lambda _1-2\mu_1)\omega
^j_1(e_1) $$ for every $j>1$. By the symmetry of $\hat\nabla K$ one
now has  $\hat\nabla_{e_1}e_1=0$. Assume now that $k>1$ and $j>1$.
Using the fact that $\hat\nabla _{e_1}e_1=0$ we obtain
$$ g((\hat\nabla _{e_j}K)(e_k,e_1), e_1)= (\lambda _1-2\mu _1)\omega
^k_1(e_j),\ \ \ \ \ \  g((\hat\nabla _{e_1}K)(e_k,e_j),e_1)=0.$$
Hence $\hat\nabla e_1=0$.  Assume now that $\hat\nabla e_1=0,...,
\hat\nabla e_{i-1}=0$. In particular, $\omega ^k_i(e_i)=0$ for every
$k<i$. By a straightforward computation one gets
$$ g((\hat\nabla _{e_j}K)(e_i,e_i),e_i)=0,\ \ \ \
g((\hat\nabla_{e_i}K)(e_j,e_i),e_i)=(\lambda_i-2\mu _i)\omega
^j_i(e_i)$$ for $j>i$. Hence $\hat\nabla _{e_i}e_i=0$.  For  $k>i$
and  any $j$ we obtain $g((\hat\nabla _{e_j}K)(e_k,e_i), e_i)
=(\lambda _i-2\mu_i)\omega ^k_i(e_j).$  In   both cases: $j>i$ and
$j<i$ one gets $g((\hat\nabla _{e_i}K)(e_k,e_j),e_i)=0$. Thus
$\hat\nabla e_i=0$ for $i=1,...,n-1$. It is now clear that
$\hat\nabla e_n=0$ as well.\koniec


\begin{thm}
Let  $(g,K)$ be a trace-free statistical structure on a manifold $M$
with symmetric $\hat\nabla K $. If the sectional $K$-curvature is
constant then either $K=0$ or $\hat R=0$ and $\hat\nabla K=0$.
\end{thm}
\proof Assume that $K\ne 0$. It means that $K_x\ne 0$ at every point
$x$ of  $M$, because the sectional curvature is constant and $K$ is
traceless. At each point of $M$ the tensor $K$ has the expression as
in Lemma \ref{ejiri} with values $\lambda _i$, $\mu _i$ given by
(\ref{wzory_na_lamba_mu}) (non-zero and constant on $M$). Moreover,
$\lambda _i- 2\mu _i\ne 0$. By Lemma \ref{smoothness_lemma} we know
that for each $p\in M$ there is a unit vector field $e_1$ in a
neighborhood of  $p$ $\, $ such that $\Phi _x$ attains a maximum
$\lambda _1$ at $(e_1)_x$  for each point $x$ of this neighborhood.
We take the orthogonal complement $\mathcal D$ to $e_1$ in the
domain of $e_1$. By Lemma  \ref{smoothness_lemma} one gets a smooth
vector field $e_2$ at which $\Phi _{|\mathcal D}$ attains a maximum
(at each point of a domain of $e_2$) and then we proceed
inductively. In this way we obtain a smooth frame field
$e_1,....,e_n$ relative to which $K$ has expression as in Lemma
{\ref{ejiri} with constant $\lambda _i$, $\mu_i$ for $i=1,..., n$.
Using now Lemma \ref{nablaKsymetryczne} completes the proof.\koniec

\begin{remark}{\rm Particular versions of the above theorem have been
given for minimal Lagrangian space forms in complex space forms, see
\cite{E} and for  affine hyperspheres  with constant sectional
curvature, see Theorem 2.2.3.18 in \cite{LSZ}.}
\end{remark}

 We shall say that a tensor $K$ of type $(1,2)$ is
non-degenerate if the mapping $X\to K_X$ is a monomorphism.

\begin{thm}
  Assume that $[K,K]=0$ on a
  statistical manifold $(M,g,K)$,  $\hat \nabla K$ is symmetric and $\hat\nabla E=0$. If $K$ is non-degenerate at each point of $M$ then
 $\hat R=0$ and $\hat \nabla K =0$ on $M$.
 \end{thm}

 \proof
 At each point $p\in M$ the tensor  $K_p$  can be expressed as in Corollary
 \ref{[K,K]=0} and all $\lambda _i$ are non-zero.
 Let $p$ be a fixed point of $M$.
 By Lemma \ref{smoothness_lemma} there is a local unit vector field $e_1$ around $p$ such that
 $\Phi$ attains its  local maximum on $U_xM$ for every $x$ from a neighborhood of $p$.
Let $\lambda_1=C(e_1,e_1,e_1)$. Take the distribution $\mathcal D$
orthogonal to $e_1$. We now take $e_2$ where $\Phi$ restricted to
$\mathcal D_p$ attains  its maximum $\lambda _2$. Again we can apply
Lemma
 \ref{smoothness_lemma}
and get a unit smooth local vector field $e_2$ in a neighborhood of
$p$ such that $\Phi_{|\mathcal D\cap U_xM}$ attains a maximum at
$e_2$ for each $x$ from this neighborhood. Continuing this process
and using the proof of Lemma \ref{ejiri} we  obtain a smooth
orthonormal local frame $e_1,...,e_n$ in a neighborhood of $p$ such
that $K(e_i,e_j)=\delta _{ij}\lambda _ie_i$ for $i,j=1,...,n$. The
functions $\lambda_i=C(e_i,e_i,e_i)$ are smooth.

We shall now use the assumption that $\hat\nabla K$ is symmetric in
order to show that $\hat \nabla e_j=0$ for all $j=1,,,,n$. For $i\ne
j$ we have
\begin{eqnarray*}
 (\hat\nabla _{e_i}K)(e_j,e_j)
  =(e_i\lambda_j)e_j+\lambda _j\sum_{l\ne j} \omega_j^l(e_i)e_l,
\end{eqnarray*}
\begin{eqnarray*}
 (\hat\nabla _{e_j}K)(e_i,e_j)= -\omega ^j_i(e_j)\lambda _je_j -\omega_j^i(e_j)\lambda_ie_i.
 \end{eqnarray*}
By comparing these equalities we get
\begin{equation}\label{1w_dowodzie_[K,K]=0}
 e_i\lambda_j=- \omega^j_i(e_j)\lambda_j,\ \ \ \ \  \ \ \ \ \omega^i_j(e_i)\lambda_j =-\omega ^i_j(e_j)\lambda_i,
\end{equation}
and
\begin{equation}
 \omega ^l_j(e_i)=0 \ \ \  for\ \  l\ne i.
\end{equation}
We now observe that $\omega ^j_i(e_j)=0$. We have $E=\lambda _1e_1
+...+ \lambda _ne_n$ and
\begin{eqnarray*}
 &&\hat\nabla_{e_i}E=(e_i\lambda _1)e_1 +... +(e_i\lambda _n)e_n\\
 &&\ \ \ \  \ \ \  \ \ \ \  \ \ \ \ \ \ \ \ \ \  +
 [\lambda _1 \omega _1^i(e_i)+...+\lambda _n\omega ^i_n(e_i)]e_i.
\end{eqnarray*}
It follows that $e_i\lambda _j=0$  for $i\ne j$. Using now
(\ref{1w_dowodzie_[K,K]=0}) we get $\omega^j_i(e_j)=0$. We have
proved that $\hat\nabla e_i=0$ for all $i=1,...,n$. In particular,
$\hat R=0$.
 Now, from the above formula for $\hat\nabla E$, we get  $e_i\lambda_i=0$. Hence all $\lambda _i$ are constant.
 It is now clear that $\hat\nabla K=0$.
 \koniec

As an immediate  consequence of the results of this paper we  have

\begin{corollary}
Let $(g,K)$ be a statistical structure on a manifold $M$  and
$\hat\nabla K=0$ on $M$. Each of the following conditions implies
that the metric $g$ is flat
\newline
1) the sectional $K$-curvature is negative
\newline
2) the sectional $K$-curvature has values in the interval
$(0,\lambda _1^2/4)$, where $\lambda _1$ is  the maximal value of
$\Phi$,
\newline
3) $[K,K]=0$ and $K$ is nondegenerate.
\end{corollary}

\end{document}